\documentclass[11pt]{article}

\usepackage{tikz}
\usepackage{subfigure}
\usepackage[english]{babel}

\usepackage[center]{caption2}
\usepackage{amsfonts,amssymb,amsmath,latexsym,amsthm}
\usepackage{multirow}
\usepackage[usenames,dvipsnames]{pstricks}
\usepackage{epsfig}
\usepackage{pst-grad} % For gradients
\usepackage{pst-plot} % For axes
\usepackage[space]{grffile} % For spaces in paths
\usepackage{etoolbox} % For spaces in paths
\makeatletter % For spaces in paths
\patchcmd\Gread@eps{\@inputcheck#1 }{\@inputcheck"#1"\relax}{}{}

\def\qedB{{\hfill\enspace\vrule height8pt depth0pt width8pt}}

\newtheorem{thm}{Theorem}%[section]

\newtheorem{lem}[thm]{Lemma}

\let\svthefootnote\thefootnote
\newcommand\blankfootnote[1]{%
	\let\thefootnote\relax\footnotetext{#1}%
	\let\thefootnote\svthefootnote%
}

\begin{document}

	\title{\bf\Large Tight bounds towards a conjecture of Gallai}
	
	\date{}
	
	\author{
		Jun Gao\thanks{Extremal Combinatorics and Probability Group, Institute for Basic Science (IBS), Daejeon, South Korea. Research supported by the Institute for Basic Science (IBS-R029-C4). Email: gj950211@gmail.com, jungao@ibs.re.kr.}~~~~~~~
		Jie Ma\thanks{School of Mathematical Sciences, University of Science and Technology of China, Hefei, Anhui 230026, China.
			Research supported by the National Key R and D Program of China 2020YFA0713100, National Natural Science Foundation of China grant 12125106,
			Innovation Program for Quantum Science and Technology 2021ZD0302904, and Anhui Initiative in Quantum Information Technologies grant AHY150200.
			Email: jiema@ustc.edu.cn.}
		%\\
		%	\medskip \\
		%	School of Mathematical Sciences\\
		%	University of Science and Technology of China\\
		%	Hefei, Anhui 230026, China.
	}

	%\blankfootnote{School of Mathematical Sciences, University of Science and Technology of China, Hefei, Anhui 230026, China.Research supported by the National Key R and D Program of China 2020YFA0713100, National Natural Science Foundation of China grant 12125106,Innovation Program for Quantum Science and Technology 2021ZD0302904, and Anhui Initiative in Quantum Information Technologies grant AHY150200.Emails: gj0211@mail.ustc.edu.cn, jiema@ustc.edu.cn.}

	\maketitle
	
	\begin{abstract}
		We prove that for $n>k\geq 3$, if $G$ is an $n$-vertex graph with chromatic number $k$ but any of its proper subgraphs has smaller chromatic number,
		then $G$ contains at most $n-k+3$ copies of a clique of size $k-1$.
		This answers a problem of Abbott and Zhou and provides a tight bound on a conjecture of Gallai.
	\end{abstract}
	
	\section{Introduction}
	A graph $G$ is called {\it $k$-critical} if its chromatic number is $k$ but any proper subgraph of $G$ has chromatic number less than $k$.
	This important notion was first introduced by G. A. Dirac in 1952 (see \cite{D52-2}) and has been extensively studied over the past decades.
	
	Throughout this paper, for a graph $G$ and a positive integer $\ell$,
	let $t_\ell(G)$ denote the number of copies of the clique $K_\ell$ on $\ell$ vertices contained in $G$.
	T. Gallai (see \cite{S87,JT95}) conjectured that every $k$-critical graph $G$ on $n$ vertices satisfies that $t_{k-1}(G)\leq n$.
	%It is clearly so when $n=k$ since complete graph of size $k$ is $k$-critical.
	This holds trivially for $k\leq 3$ (note that the only $3$-critical graphs are odd cycles).
	Using an elegant argument of linear algebra, Stiebitz~\cite{S87} confirmed this for $4$-critical graphs $G$ by showing $t_3(G)\leq n$.
	On the other hand, he~\cite{S87} proved that for any $k\geq 4$, there exist some constant $c_k>0$ and arbitrarily large $k$-critical graphs $G$ on $n$ vertices
	such that $t_\ell(G)\ge c_k n^\ell$ holds for each $\ell\in \{2,3,...,k-2\}$.
	Koester~\cite{K91} provided an improvement for $4$-critical planar graphs $G$ by showing that if $G$ has $n\geq 6$ vertices, then $t_3(G)\leq n-1$.
	The cases $k\ge 5$ of Gallai's conjecture were resolved completely by Abbott and Zhou in \cite{AZ92},
	who extended Stiebitz's arguments and proved that any $k$-critical graph $G$ on $n$ vertices has $t_{k-1}(G)\leq n$ with equality only if $n=k$ and $G\cong K_n$.
	They \cite{AZ92} also showed that for any $4$-critical graph $G$ on $n$ vertices,
	if $G$ is not an odd wheel\footnote{An {\it odd wheel} is obtained from an odd cycle $C$ by adding a new vertex $x$ and joining $x$ to every vertex of $C$. Note that an odd wheel on $n\geq 6$ vertices is a $4$-critical planar graph and has exactly $n-1$ triangles.}, then $t_3(G)\leq n-2$.
	For integers $\ell, d\geq 2$, let $W(\ell,d)$ denote the graph obtained from a disjoint union of a clique $K_d$ on $d$ vertices and a cycle $C_\ell$ of length $\ell$ by joining each vertex of $K_d$ to each vertex of $C_\ell$.
	Observe that if $n-k+3$ is odd, then $W(n-k+3,k-3)$ is an $n$-vertex $k$-critical graph with exactly $n-k+3$ copies of $K_{k-1}$.
	Abbott and Zhou~\cite{AZ92} posed the following problem, which was stated as a conjecture in K\'ezdy and Snevily~\cite{KS97}.
	
	\medskip
	
	{\bf \noindent Conjecture.} (Abbott and Zhou \cite{AZ92})
	Let $G$ be an $n$-vertex $k$-critical graph with $n>k\geq 4$. Then $t_{k-1}(G)\leq n-k+3$.
	
	\medskip
	
	This (if true) would be tight for infinitely many integers $n$ as indicated by the above graph $W(n-k+3,k-3)$.
	The aforementioned result of Abbott and Zhou~\cite{AZ92} on $4$-critical graphs implies the case $k=4$,
	and the cases $k\leq 7$ were confirmed by Su~\cite{S95,Su95}.
	The current best bound for the general case was obtained by Kezdy and Snevily~\cite{KS97} as follows.
	
	\begin{thm}[K\'ezdy and Snevily~\cite{KS97}]
		Let $G$ be an $n$-vertex $k$-critical graph with $n>k\geq 4$. Then $t_{k-1}(G)< n-3k/5+2$.
	\end{thm}
	
	The proof of this theorem uses linear algebra as well as some careful analysis from structural graph theory.
	We mention that the above problems and results are discussed in detail in Section 5.9 of the book of Jensen and Toft~\cite{AZ92} (see its page 103).
	
	In this paper, we confirm the conjecture of Abbott and Zhou by proving the following.
	\begin{thm}\label{main}
		Let $n>k\geq 4$. Any $n$-vertex $k$-critical graph $G$ has $t_{k-1}(G)\leq n-k+3$.
	\end{thm}
	
	Our proof uses linear algebra arguments, which originate from Stiebitz~\cite{S87} and appear in the subsequent works \cite{AZ92,KS97}.
	We would like to emphasize that the core part of our proof is different from \cite{KS97}, which we will elaborate in Section~2.
	
	\section{The proof}
	To present the proof of Theorem~\ref{main}, we will first need to introduce some notation and several existing results.
	Let $G$ be an $n$-vertex graph with vertex set $V(G)=\{v_1,v_2,..., v_n\}$.
	For a subset $S\subseteq V(G)$, we define its {\it incidence vector} to be a 0-1 vector $\vec{u}_S= (u_1,u_2,...,u_n)$,
	where $u_i=1$ if $v_i\in S$ and $u_i=0$ otherwise.
	The first lemma we need is given by Stiebitz~\cite{S87}, which reveals the special role of the graph $W(\ell,k-3)$ in $k$-critical graphs.

	\begin{lem}[Stiebitz~\cite{S87}]\label{lem:Sti}
		Let $k\geq 4$. If $G$ is a $k$-critical graph containing some $W(\ell, k-3)$ as a subgraph, then $G\cong W(\ell, k-3)$ and $\ell$ is an odd integer.
	\end{lem}

	The following lemma is a direct consequence of a result of Abbott and Zhou~\cite{AZ92}.
	
	\begin{lem}[Abbott and Zhou~\cite{AZ92}, see its Lemma~2]\label{lem:AZ}
		Let $k\geq 4$ and $G$ be a $k$-critical graph that does not contain any $W(\ell, k-3)$ as a subgraph.
		Let $\vec{x}_1,\vec{x}_2,...,\vec{x}_r$ be incidence vectors of all cliques $K_{k-1}$ in $G$.
		Then $\vec{x}_1,\vec{x}_2,...,\vec{x}_r$ are linearly independent over $GF(2)$.
	\end{lem}
	
	%\begin{lem}[Abbott and Zhou~\cite{AZ92}, see Lemma~2]\label{linear}
	%		Let $k\ge 3$ and let $G$ be a graph with no isolated vertices. Suppose that
	%		\begin{enumerate}
		%			\item each edge of $G$ is an edge of some $K_{k-1}$,
		%			\item G does not contain $W(\ell,k-3)$, for any $\ell\ge3$,
		%			\item for each edge $(a,b)$ of $G$, the graph $G-(a,b)$ has a $(k-1)$-coloring in which $a$ and $b$ have the same color, then :
		%		\end{enumerate}	
	
	%	(A) there is a vertex $v$ of $G$ such that the number of $(k-1)$-cliques containing $v$ is odd,
	
	%	(B) the number $t$ of $(k-1)$-cliques of $G$ satisfies $t\le n-1$, where $n$ is the order of $G$.
	%	\end{lem}

%\begin{lem}[Stiebitz~\cite{S87}]\label{S:lem}
%If $G$ is a $k$-critical graph ($k\ge4$) containing a $W(\ell,k-3)$ as a subgraph, then $G \cong W(\ell,k-3)$ and $\ell$ is an odd integer.
%\end{lem}

The following nice lemma relates the total number of cliques in a $k$-critical graph to the number of cliques containing any fixed edge.
The cases $d\in \{0,1\}$ were first obtained by Su~\cite{Su95} and the general case was later proved by Kezdy and Snevily~\cite{KS97}.
We shall mention that our proof will only need the case $d=0$.

\begin{lem}[Kezdy and Snevily~\cite{KS97}]\label{KS:lem}
	Let $G$ be an $n$-vertex $k$-critical graph.
	If there is an edge in $G$ that is contained in exactly $d$ copies of $K_{k-1}$, then $t_{k-1}(G)\leq n-(k-2-d)$.
\end{lem}

%Let $G$ be a graph. For $S\subseteq V(G)$, let $G[S]$ be the subgraph of $G$ induced on the vertex set $S$.

\medskip

We are ready to present the proof of our result Theorem~\ref{main}.

\medskip

{\noindent \bf Proof of Theorem~\ref{main}.}
Let $n>k\geq 4$ be integers and let $G$ be any $k$-critical graph on $n$ vertices.
We aim to show that $t_{k-1}(G)\leq n-k+3$.

If $G$ contains some $W(\ell, k-3)$, then by Lemma~\ref{lem:Sti},
$G\cong W(\ell, k-3)$ and $\ell=n-k+3\geq 4$ is odd, from which the desired conclusion $t_{k-1}(G)=n-k+3$ holds.
Hence we may assume that there is no copy of $W(\ell, k-3)$ in $G$.
In particular there is no $K_k$ in $G$ and $G\not\cong K_k$,
so by the result of Abbott and Zhou \cite{AZ92}, we have $t_{k-1}(G)\leq n-1$.
%By Lemma~\ref{S:lem}, we assume $G$ does not contain a $W(\ell,k-3)$ as a subgraph for any $\ell \ge 3$.
For any $x\in V(G)$, let $t_{k-1}(x,G)$ denote the number of cliques $K_{k-1}$ in $G$ containing $x$.
Then we have $$\sum_{x\in V(G)} t_{k-1}(x,G) = (k-1)\cdot t_{k-1}(G) \leq (k-1)(n-1).$$
Let $u$ be the vertex minimizing $t_{k-1}(u,G)$ among all vertices in $G$.
By the above inequality, we see that $t_{k-1}(u,G)\leq k-2$.

Suppose that the neighborhood $N(u)$ of the vertex $u$ induces a complete subgraph of $G$.
In this case, as $G$ does not contain any copy of $K_k$, we see $|N(u)|\leq k-2$.
This is a contradiction, as the minimum degree of a $k$-critical graph is at least $k-1$.

Therefore, there exist two vertices $v,x\in N(u)$ such that $v,x$ are not adjacent in $G$.
For any edge $e\in E(G)$, we denote $t_{k-1}(e,G)$ to be the number of copies of $K_{k-1}$ in $G$ that contain $e$.
We may assume that $t_{k-1}(e,G)\geq 1$, i.e., any edge $e$ is contained in at least one copy of $K_{k-1}$ (as otherwise Lemma~\ref {KS:lem} implies that $t_{k-1}(G) \leq n-k+2$).

Since $v,x$ are not adjacent, the set of all cliques $K_{k-1}$ containing $uv$ is disjoint from the set of all cliques $K_{k-1}$ containing $ux$.
So we have $$t_{k-1}(uv,G)+t_{k-1}(ux,G)\leq t_{k-1}(u,G)\leq k-2.$$
Because $t_{k-1}(ux,G)\geq 1$, we see that
\begin{equation}\label{equ:t(uv)}
	t_{k-1}(uv,G)\leq k-3.
\end{equation}

\noindent\textbf{Claim 1.} There exists a clique $A=\{a_1,a_2,...,a_{k-1}\}$ of size $k-1$ such that $x\in A$ and $A\cap \{u,v\}= \emptyset$.

\smallskip

\begin{proof}
	First we have by the minimality of $t_{k-1}(u,G)$ that $t_{k-1}(x,G) \ge t_{k-1}(u,G)$.
	There also exists a clique $K_{k-1}$ containing $uv$, that, in particular, contains $u$ but not $x$.
	These two facts together indicate that
	there exists a clique $A$ of size $k-1$, that contains $x$ but not $u$.
	As $xv\notin E(G)$, this clique $A$ cannot contain $v$ as well, proving the claim.
\end{proof}

Because $G$ is $k$-critical, the subgraph $G-uv$ is $(k-1)$-colorable and thus
there exists a proper coloring $\phi:V(G)\to \{1,2,..., k-1\}$ of $G-uv$
such that $u$ and $v$ are assigned the same color, say the color $k-1$.
We now prove the following claim.

\medskip

\noindent\textbf{Claim 2.} There exists a color $c\in \{1,2,...,k-2\}$ such that
every clique $K_{k-1}$ in $G$ contains a vertex that is colored by $c$ under $\phi$.
We may assume $c=k-2$.

\smallskip

\begin{proof} To see this, we first note that each of the cliques $K_{k-1}$ in $G$ not containing the edge $uv$
	must use all colors in $\{1,2,..., k-1\}$ under $\phi$.
	For any clique $K_{k-1}$ in $G$ containing the edge $uv$,
	it uses exactly $k-3$ colors in $\{1,2,...,k-2\}$ under $\phi$,
	for which one color needs to be removed from the list $\{1,2,...,k-2\}$ for this claim.
	By \eqref{equ:t(uv)}, there are at most $k-3$ such cliques $K_{k-1}$,
	which together will remove at most $k-3$ colors from the list $\{1,2,...,k-2\}$ for this claim.
	This leaves at least one color $c\in \{1,2,...,k-2\}$ such that every clique in $G$ witnesses the color $c$ under $\phi$.
\end{proof}

For each $1\leq i\leq k-1$, let $$C_i=\{x\in V(G): \phi(x)=i\}.$$
For the clique $A=\{a_1,a_2,...,a_{k-1}\}$ from Claim~1, we may assume that $a_i\in C_i$.
Let us recall the properties of $A$ and it will be crucial for us to notice that
\begin{equation}\label{equ:a(k-1)}
	a_{k-1}\in C_{k-1}\backslash \{u,v\}.
\end{equation}
Let $r= t_{k-1}(G)$ and let $T_1,T_2,\cdots,T_r$ be all cliques $K_{k-1}$ in $G$.
%$V(G)=\{v_1,v_2,\cdots,v_n\}$
For each $1\leq i\leq r$, we use $\vec{x}_i$ to denote the incidence vector of $T_i$,
and for each $1\leq j\leq k-3$, we use $\vec{y}_j$ to denote the incidence vector of the single-vertex set $\{a_j\}$.\footnote{Note that here we only use $k-3$ incidence vectors from $A$ to form $\vec{y}_j$'s. In total, there are $k-1$ elements of $A$ that correspond to $k-1$ colors. We have two special colors $k-1$ and $k-2$ set aside after Claim 1 and Claim 2, respectively, which leaves $k-3$ colors.}

The rest of the proof will be devoted to show the statement that
$$\mbox{the vectors $\vec{x}_1,\vec{x}_2,...,\vec{x}_r,\vec{y}_1,\vec{y}_2,...,\vec{y}_{k-3}$ are linearly independent over $GF(2)$.}$$
Note that all these vectors are defined in an $n$-dimensional linear space over $GF(2)$.
If this statement is proved to be true, then we have $r+(k-3)\leq n$,
from which the conclusion of Theorem~\ref{main} that $t_{k-1}(G)=r\leq n-k+3$ holds.

Suppose for a contradiction that there exist $\vec{x}_{q_1},...,\vec{x}_{q_s},\vec{y}_{p_1},...,\vec{y}_{p_m}$ such that
\begin{equation}\label{equ:vecs}
	\vec{x}_{q_1}+\vec{x}_{q_2}+\cdots +\vec{x}_{q_s}+\vec{y}_{p_1}+\vec{y}_{p_2}+\cdots+\vec{y}_{p_m}=\vec{0}
\end{equation}
over $GF(2)$, where $1\leq q_i \leq r$ and $1\leq p_j\leq k-3$ for all possible $1\leq i\leq s$ and $1\leq j\leq m$.
We may assume $q_i =i$ and $p_j =j$ for all $i$ and $j$.
By Lemma~\ref{lem:AZ}, $\vec{x}_1,\vec{x}_2,...,\vec{x}_r$ are linearly independent over $GF(2)$,
so $m\geq 1$. Since the $\vec{y_i}$ are independent as well, we have $s\geq 1$.

Let $\mathcal{G} = \{T_1,T_2,\cdots,T_s\}$.
For any vertex $w\in V(G)$ and any pair $e\in \binom{V(G)}{2}$, let $t(w,\mathcal{G})$ and $t(e,\mathcal{G})$ denote the number of cliques in $\mathcal{G}$ containing $w$ and $e$, respectively.\footnote{If $e\notin E(G)$, then it is evident that we have $t(e,\mathcal{G})=0$.}
We observe from \eqref{equ:vecs} that for any vertex $w\in V(G)$,
\begin{equation}\label{equ:odd-t(w)}
	\mbox{ $t(w,\mathcal{G})$ is odd if and only if $w\in \{a_1,a_2,...,a_m\}$.}
\end{equation}
As $1\leq m\leq k-3$, we get that $\{a_1,a_2,...,a_m\}\cap (C_{k-2}\cup C_{k-1})=\emptyset$,
showing that $t(w,\mathcal{G})$ for all $w\in C_{k-2}\cup C_{k-1}$ are even (in particular, $t(a_{k-1},\mathcal{G})$ is even).
By Claim 2, every clique in $\mathcal{G}$ contains exactly one vertex in $C_{k-2}$,
so we derive that
\begin{equation}\label{equ:parity-s}
	\mbox{ $|\mathcal{G}|=|\{(w,T_j): w\in C_{k-2}\cap T_j \mbox{ and } T_j\in \mathcal{G}\}|=\sum_{w\in C_{k-2}} t(w,\mathcal{G})$ is even. }
\end{equation}

We have seen from \eqref{equ:odd-t(w)} that $t(a_{k-1},\mathcal{G})$ is even.
To reach the final contradiction,
we want to estimate the parity of $t(a_{k-1},\mathcal{G})$ using a different approach, i.e., by looking at the contributions of all edges between $a_{k-1}$ and $C_1$.
This will be done in the coming claim.

\medskip

\noindent\textbf{Claim 3.} $t(a_1a_{k-1},\mathcal{G})$ is odd, and for any $w\in C_1\backslash \{a_1\}$, $t(wa_{k-1},\mathcal{G})$ is even.	

\smallskip

\begin{proof}
	Let $w\in C_1$ be any vertex.
	If $wa_{k-1}\notin E(G)$, then it is clear that $w\neq a_1$ and $t(wa_{k-1},\mathcal{G})=0$ that is even.
	So from now on we may assume $wa_{k-1}\in G$.
	Then there exists a proper coloring $\chi_w: V(G)\to \{1,2,..., k-1\}$ of $G-wa_{k-1}$
	such that $w$ and $a_{k-1}$ are assigned the same color, say the color $k-1$.
	It is easy to see that every clique $K_{k-1}$ containing the edge $wa_{k-1}$ has exactly two vertices
	(i.e., $w$ and $a_{k-1}$) with the color $k-1$ under $\chi_w$,
	while every other clique $K_{k-1}$ not containing $wa_{k-1}$ has exactly one
	vertex with the color $k-1$ under $\chi_w$; let us call this property $(\star)$.

	Suppose that $w\in C_1\backslash \{a_1\}$.
	For any vertex $z$ with $\chi_w(z) = k-1$,
	we have $z=a_{k-1}$, or $z=w\in C_1\backslash \{a_1\}$, or $z$ is not adjacent to $a_{k-1}$.
	However in any case, such $z$ is not contained in $\{a_1,a_2,...,a_m\}$.
	By \eqref{equ:odd-t(w)}, we see that $t(z,\mathcal{G})$ is even for any vertex $z$ with $\chi_w(z) = k-1$.
	Let $\Lambda$ denote the number of pairs $(z,T_j)$ satisfying $z\in T_j\in \mathcal{G}$ and $\chi_w(z) = k-1$.
	By the above fact and the property $(\star)$, we get that
	$$|\mathcal{G}|+t(wa_{k-1},\mathcal{G})=\Lambda=\sum_{z:~\chi_w(z) = k-1} t(z,\mathcal{G}) \mbox{ is even.}$$
	As $|\mathcal{G}|$ is even (i.e., \eqref{equ:parity-s}), we then derive that $t(wa_{k-1},\mathcal{G})$ is even for any $w\in C_1\backslash \{a_1\}$.

	It remains to consider $w=a_1$.
	By similar analysis, for any vertex $z\neq w$ with $\chi_w(z) = k-1$,
	we get that $t(z,\mathcal{G})$ is even.
	This together with the fact that $t(w,\mathcal{G})$ is odd imply that
	$$|\mathcal{G}|+t(wa_{k-1},\mathcal{G}) = \sum_{z:~\chi_w(z) =k-1} t(z,\mathcal{G}) \mbox{ is odd.}$$
	Thus $t(a_1a_{k-1},\mathcal{G})=t(wa_{k-1},\mathcal{G})$ is odd. This finishes the proof of Claim~3.
\end{proof}

By fact \eqref{equ:a(k-1)},
every clique $K$ of size $k-1$ in $G$ containing $a_{k-1}$ does not contain the edge $uv$,
so we have $|K\cap C_i|=1$ for each $i\in \{1,2,...,k-1\}$.
This shows that $$t(a_{k-1},\mathcal{G}) = \sum_{w\in C_{1}} t(wa_{k-1},\mathcal{G}),$$
which together with Claim~3 imply that
$t(a_{k-1},\mathcal{G})$ is odd.
But this is a contradiction to \eqref{equ:odd-t(w)} as it oppositely says that $t(a_{k-1},\mathcal{G})$ is even.
The proof of Theorem~\ref{main} now is complete.\qedB

\medskip

It would be very interesting to determine all $n$-vertex $k$-critical graphs $G$ with $t_{k-1}(G)=n-k+3$ for $n>k\geq 4$.
For the case $k=4$, it is known from the result of Abbott and Zhou \cite{AZ92} that such $4$-critical graphs can only be odd wheels.
For $k\geq 5$, it seems to be challenging to say something about the structure of these $k$-critical graphs from the proof presented here.
We tend to believe that in case $n-k+3$ is odd, the graph $G=W(n-k+3, k-3)$ is the only extremal graph for Theorem~\ref{main} satisfying that $t_{k-1}(G)=n-k+3$.

There also is a related conjecture proposed by Su~\cite{Su95}, which states that any $k$-critical graph of order $n>k$ has an edge that is contained in at most one clique $K_{k-1}$ on $k-1$ vertices.
Su proved that this conjecture would imply Theorem~\ref{main}, and this proof was extended by K\'ezdy and Snevily~\cite{KS97} to Lemma~\ref{KS:lem}.
The cases $4\le k\le 7$ were verified by Su~\cite{Su95}.
It is interesting to have an alternative proof of Theorem~\ref{main} via this conjecture.

\bigskip

\noindent{\bf Acknowledgement.}  The authors would like to thank two referees and Stijn Cambie for their careful reading and for many valuable suggestions.

\end{document}